\documentclass[14]{amsart}
\usepackage{graphicx}
\usepackage{amssymb}
\usepackage{epstopdf}
\usepackage{nicefrac}
\usepackage{enumerate}
\usepackage{hyperref}
\usepackage{color}
\usepackage{pst-all}
\usepackage{tikz}
\usepackage{pgfplots}
\usepackage{float}
\usepackage{pstricks, pst-plot, pst-3dplot,pst-solides3d}

\usepackage [autostyle, english = american]{csquotes}
\MakeOuterQuote{"}

\graphicspath{{./Figures/}}

\theoremstyle{plain} 
\newtheorem{theorem}{Theorem}[section]
\newtheorem{lemma}[theorem]{Lemma}

\newtheorem*{theorem*}{Theorem}
\newtheorem*{lemma*}{Lemma}
\newtheorem*{exercise*}{Exercise}

\theoremstyle{definition} 

\newtheorem*{definition*}{Definition}
\newtheorem*{notation*}{Notation}
\newtheorem*{example*}{Example}
\newtheorem*{question*}{Question} 

\newtheorem*{answer*}{Answer}

\title{ Higher order  Petri Loci.}
\author{Montserrat Teixidor i Bigas}
\address{Mathematics Department, Tufts University\\%
	177 College Ave\\%
	Medford, MA 02155\\%
	USA\\%
	ORCiD 0000-0003-3747-3330}
	\email{mteixido@tufts.edu}

\begin{document} 
\begin{abstract} Denote by  ${\mathcal P}_{g,d}^{r,k}$ the subset of ${\mathcal M}_g$ consisting of curves that have a linear series of degree $d$ and dimension $r$ for which the Petri map has kernel of dimension at least $k$.
We show the existence of codimension $k$ components of  ${\mathcal P}_{g,d}^{r,k}$.
\end{abstract}
\maketitle
\begin{center}

\end{center}
\vskip .5cm

Let $L$ be a line bundle on a curve $C$ of genus $g$ and   $W$  a  subspace of the space of sections of $L$, $W\subseteq H^0(C, L)$.
The Petri map is the natural cup-product 
\[ \mu_0:W\otimes H^0(C, K\otimes L^{*})\to H^0(C, K\otimes L\otimes L^{*})\simeq  H^0(C, K) \]
For a generic curve $C$ in  the moduli space ${\mathcal M}_g$ of curves of genus $g$, the Petri map is known to be injective for every choice of $L$ and  $W\subseteq H^0(C, L)$.

Denote by  ${\mathcal P}_{g,d}^{r,k}$ the necessarily proper subset of ${\mathcal M}_g$ consisting of curves that have a linear series of degree $d$ and dimension $r$ for which the Petri map has kernel of dimension at least $k$.
These loci have been studied when $k=1$, that is when the Petri map fails to be injective.  
When  ${\mathcal P}_{g,d}^{r,1}$  are non-empty for a given $r,d$, they are expected to be divisors. 
In fact, it is know that this is the case  for positive values of the Brill-Noether number $\rho\ge 0 $,
 if   ${\mathcal P}_{g,d}^{r,1}$ is not contained in  ${\mathcal P}_{g,d}^{r+1,1}$ (see   \cite{BS}).
 The proof of this result is based on the natural determinantal structure of   ${\mathcal P}_{g,d}^{r,1}$:
 Let us denote by  ${\mathcal G}_{g,d}^r$ the scheme parameterizing linear series of degree $d$ and dimension $r$,
 which can be defined locally in the neighborhood of a curve and line bundle.
 From Brill-Noether's Theorem, for generic $C$ and  in a neighborhood of a $(C,L,W)$,  ${\mathcal G}_{g,d}^r$ has dimension $\dim {\mathcal M}_g+\rho$.
 The Petri locus can be described by the vanishing of the higher order minors of the matrix of $\mu _0$.
 Therefore, the expected codimension  of the set of triples $(C, L, W)$ where the Petri map fails to be injective is $g-[(r+1)(g-d+r)-1]=\rho +1$.
  If the projection to ${\mathcal M}_g$ is transversal, the codimension  of ${\mathcal P}_{g,d}^{r,1}$ in ${\mathcal M}_g$ is 1.
   For small values of the genus, ad hoc methods allow to show that every component of   ${\mathcal P}_{g,d}^{r,1}$ has codimension one (see \cite{C}, \cite{LC1}).

 It had been conjectured that any component of  ${\mathcal P}_{g,d}^{r,1}$  should be a divisor  in ${\mathcal M}_g$.
  This was disproved in \cite{LC2} which displays  a  codimension two locus   ${\mathcal P}_{18,16}^{2,1}\subseteq {\mathcal M}_{18}$ 
  and is   not contained in any codimension one locus of ${\mathcal P}_{18,d}^{r,1}\subseteq {\mathcal M}_{18}$ (even when allowing for other values of $d, r$).
  In this case  ${\mathcal P}_{18,16}^{2,1}= {\mathcal P}_{18,16}^{2,2}$, that is, the kernel of the Petri map has dimension 2 which is also $-\rho$.
 We expect this to be the rule rather than the exception: when $\rho$ is negative,  we denote with ${\mathcal M}^r_{g,d}$  the proper subset of  ${\mathcal M}_g$ 
  consisting of curves having a linear series of degree $d$ and rank $r$.
 As the  Brill-Noether number can be written as $ \rho=\dim H^0(C, K)-\dim(H^0(C,L)\otimes H^0(C, KL^{*}))$
  if a curve has a line bundle for a value of $r,d$ with negative $\rho$, then the Petri map is necessarily non-injective,
or more precisely,  ${\mathcal M}^r_{g,d}\subseteq {\mathcal P}_{g,d}^{r,-\rho}$.
This of course does not discard the possibility that  ${\mathcal M}^r_{g,d}$ could be contained in some divisorial component of $ {\mathcal P}_{g,d'}^{r',1}$ for some values of $d',r'$.

We showed in  \cite{distBN}, Theorem 2.1 that when the Brill Noether number $\rho$ is negative but $\rho\ge r+1-g$, 
then there exists a subset of  codimension $-\rho$ of ${\mathcal M}_g$ of curves possessing a $g^r_d$.
Moreover the Petri map on the generic curve of the locus is onto, so that the kernel of the Petri map has  dimension precisely $-\rho$.
In particular, the codimension of the locus possessing these linear series coincides with the dimension of the kernel of the Petri map.

Something similar holds for curves with theta characteristics with many sections. 
We say that a curve $C$ has a theta characteristic of dimension $r$ if it has a line bundle $L$ such that $L^2=K_C$ with $\dim H^0(C,L)=r+1$.
A theta characteristic of dimension $r$ is expected to have a kernel of the Petri map of dimension ${r+1\choose 2}$ spanned by elements of the form $s_i\otimes s_j-s_j\otimes s_i$ 
for $s_i, s_j$ distinct sections of $L$.
The loci of curves with a theta characteristic of dimension at least $r$ can be represented as a degeneracy locus of a symmetric quadratic form (see \cite{H}).
As such, its expected codimension is ${r+1\choose 2}$.
Equality between the codimension  in  ${\mathcal M}_{g}$ of the locus of curves with theta-characteristics and the dimension of the kernel of the Petri map holds in many cases 
(see \cite{semican},\cite{Benzo}, \cite{HT}).

Both Brill-Noether loci for negative $\rho$ and loci of curves with $r$-dimensional theta characteristics provide  two series of examples of ${\mathcal P}_{g,d}^{r,k}$  of codimension $k$ in ${\mathcal M}_g$.
The natural  description of  ${\mathcal P}_{g,d}^{r,k}$ as a determinantal variety in  ${\mathcal G}_{g,d}^r$  for $k>1$ 
would suggest a codimension much larger than $k$ in   ${\mathcal M}_g$. There  is  though another point advocating for the codimension $k$ option:
Recall that the orthogonal to the image of the Petri map is identified with the tangent space to $W^r_d$, that is,
 it parameterizes deformations of a line bundle to which the space of sections can be extended.
There is a second Petri map $\mu_1$ from the kernel of $\mu_0$ to $H^1(C, T_C)$.
The orthogonal to the image parameterizes deformations of the curve to which the kernel of the Petri map  extends (see section \ref{secdeff}). 
If $\mu_1$ is injective and the codimension of ${\mathcal P}_{g,d}^{r,k}$ in  ${\mathcal M}_g$  is $k$, we obtain a smoothness result for  ${\mathcal P}_{g,d}^{r,k}$.
We will show the existence of smooth components of codimension $k$ of   ${\mathcal P}_{g,d}^{r,k}$ that extend beyond the ${\mathcal M}^r_{g,d}$ for negative values of $\rho$ (see Theorem \ref{ggeneral}).
We use  constructive procedures and degeneration methods similar to the ones used to show the non-emptiness of  ${\mathcal P}_{g,d}^{r,1}$ in \cite{F1}, \cite{CT}  \cite{F2}.

On the other hand, there are loci  ${\mathcal P}_{g,d}^{r,k}$  of codimension smaller than $k$.
Take for example  ${\mathcal P}_{g,4}^{2,k} , 1\le k\le 2g-6$. Then  ${\mathcal P}_{g,4}^{2,k} ={\mathcal M}^1_{g,2}$ which has codimension $g-2$ in ${\mathcal M}_{g}$.

\section{Deformation of curves and line bundle }\label{secdeff}
Let $C$ be a curve, $L$ a line bundle on $C$.
Denote by $T_C$  the sheaf of derivations.
A  first order differential operator $\sigma$  acting on $L$ satisfies $\sigma (fx)=\nu(f)s+f \sigma(s)$ for some derivation $\nu$. 
 Denote  by $\Sigma_L$ the sheaf of  first order differential operators acting on $L$.
There exists an exact sequence 
\[ 0\to {\mathcal O}_C\to\Sigma_L\to T_C\to 0   \]
where $i: {\mathcal O}_C\to\Sigma_L$ sends a function $f$ to the differential operator ``multiplication with $f$'' 
and  $\pi: \Sigma_L\to T_C$ assigns to a differential operator $\sigma$ the corresponding derivation $\nu$.

The space $H^1(C, O_C)=H^1(C, L\otimes L^*)$ parameterizes first order deformations of the line bundle $L$,
the space $H^1(C, T_C)$ parameterizes first order deformations of the curve while the space $H^1(C, \Sigma_L)$ parameterizes first order deformations of the pair consisting of the curve and the line bundle.
The exact sequence above gives rise to an exact sequence 
\[ 0\to H^1(C, {\mathcal O}_C)\to  H^1(C, \Sigma_L) \to  H^1(C, T_C)\to 0   \]
With the identifications of these spaces with spaces of deformations, the first map sends a deformation of a line bundle to a deformation of the pair  when the curve stays fixed.
The second map is the  forgetful map  that ignores the deformation of the bundle and  remembers only the deformation of the curve.

Consider  a  subspace $W$ of the space of sections of $L$, $W\subseteq H^0(C, L)$.
The deformations of $L$ have extensions of the sections in $W$  if they are  in the orthogonal to the image of the Petri map
\[ \mu_0:W\otimes H^0(C, K\otimes L^{*})\to H^0(C, K\otimes L\otimes L^{*})\]
which by duality becomes 
\[ H^1(C,  L\otimes L^{*})\to Hom(W,  H^1(C, L)) \]
Similarly, the deformations of the pair $(C, L)$ that preserve a certain subspace $W\subseteq H^0(C, L)$ are those in the kernel of the cup-product map
 $P_1^*:H^1(C, \Sigma_L)\to Hom(W, H^1(C,L))$.

There is  a commutative diagram 
\[ \begin{matrix}0&\to&H^1(C, {\mathcal O}_C)&\to&H^1(C, \Sigma_L)&\to&H^1(C, T_C)&\to& 0 \\  & &\downarrow&  &\downarrow&  &\downarrow& &\\
 & & Hom(W, H^1(C,L))& =& Hom(W, H^1(C,L))&\to&(Ker \mu_0)^*& & \end{matrix}  \]
 The dual of the last vertical map often called the Wahl map was actually introduced by Arbarello and Cornalba \cite{AC}:
 \[\begin{matrix} \mu_1:& Ker \mu_0& \to &H^0(C, 2K)\\ &\sum w_i\otimes v_i&\to & \sum w_i d(v_i)\end{matrix}\]
The orthogonal to the image of $\mu_1$ parameterizes deformations of the curve to which the kernel of the Petri map  extends. 
It can therefore be interpreted as  the tangent space to ${\mathcal P}_{g,d}^{r,k}$.
If the codimension of ${\mathcal P}_{g,d}^{r,k}$ in ${\mathcal M}_g$ coincides with the dimension of the image of $\mu_1$, 
then ${\mathcal P}_{g,d}^{r,k}$ is smooth of this codimension in ${\mathcal M}_g$.

The map $\mu_1$ can be described explicitly when $\dim W=2$, and the  $g^1_k$ spanned by the sections $s_1,s_2$ has no fixed points (see Lemma 2.12 of \cite{semican}):
From the base-point-free pencil trick, the kernel of the Petri map can  then be identified to $H^0(C, K-2L)$ as follows
\[ \begin{matrix}H^0(C, K-2L)&\to& H^0(C, L)\otimes H^0(C, K-L)\\ t&\to &s_1\otimes s_2t-s_2\otimes s_1t\end{matrix}\]
The action of $(\sigma_{ij})\in H^1(C, \Sigma_L)$ is given by 
\[ \begin{matrix}H^1(C, \Sigma_L)&\to &Hom(H^0(C,L), H^1(C,L))\\ (\sigma_{ij})&\to&(a_i)\to (\sigma_{ij}(a_i))\end{matrix}\]
As $s_1, s_2$ are sections of the same line bundle $L$, there is a function $f$ such tht $s_1=fs_2$.
Then
\[ \sigma(s_1)ts_2-\sigma(s_2)ts_1=\sigma(fs_2)ts_1-\sigma(s_1)tfs_1=\nu(f)ts_1^2.\]
Note that the condition $s_1=fs_2$ means that $f$ gives the map $C\to {\mathbb P}^1$.
Then $\nu(f)$ vanishes on the ramification divisor of $f$.
The map 
\[ \begin{matrix}H^1(C, T_C)&\to &(Ker P)^*= ( H^0(C, K-2L) )^*\\ (\nu_{ij})&\to&t\to \nu_{ij}(f) ts_1^2\end{matrix}\]
gives rise to sections vanishing on $R$.

In particular we obtain the following:  
\begin{lemma}\label{tnkgonal} Let $k, g$ be integers such that $2d-g-2<0$ or equivalently, the generic curve of genus $g$ is not $k$-gonal.
 Let $C$ be a $k$-gonal curve and assume that the  $g^1_k$ is complete and without fixed points and let  $R$ be its ramification divisor.
 The tangent space to the $k$-gonal locus at C can be identified with the orthogonal to $H^0(C,2K-R)$ where we consider $H^0(C,2K-R)\subseteq H^0(C,2K)\cong H^1(C,T_C)^*$.
\end{lemma}

\section{Limit Linear series in a chain of elliptic curves.}

Given elliptic curves $C_1,\dots, C_g$  and points $P_i, Q_i\in C_i$,  form a nodal curve of genus $g$ by identifying $Q_i$ to $P_{i+1}, i=1,\dots, g-1$.
Assume that the curves and chosen points are generic except for a few of the curves  $C_{i_1},\dots,C_{i_e}$  for which there exists an integer $l_j$ such that $P_{i_j}-Q_{i_j}$ is a point on the elliptic curve with $l_j$ torsion.
In \cite{distBN}, Prop 1.4 we gave a description of the smoothable limit linear series on such chains:

The components of the space of the smoothable limit linear series on the chain correspond to the fillings of an $(r+1)\times (g-d+r)$ rectangle with indices $1,\dots, g$
 such that both rows and columns are in increasing order,
 if an index appears more than once say on spots $(a_1,b_1)$ and $(a_2, b_2)$,  then  the grid distance between the two $(a_2-a_1)+(b_1-b_2)$ is divisible by $l_i$ 
 and no index could be replaced by another still satisfying the conditions to allow for fewer repetitions. 
 We called this an admissible filling of the rectangle.
  To simplify notations, we will often write $\alpha=r+1, \beta=g-d+r$.

Given a limit linear series on the chain corresponding to line bundles $L_i$ and  spaces of sections $V_i$, a section of the limit series is a collection of sections   $(s^i) \ i=1,\dots, g, s^i\in Vi$ such that 
 $ord_{Q^i}s^i+ord_{P^{i+1}}s^{i+1}\ge d$. 
With the choice of linear series corresponding to the filling of the $\alpha\times \beta$ rectangle, there are   sections $s_1,\dots, s_{\alpha}$ of the linear series 
and sections $t_1,\dots, t_{\beta}$  of the Serre dual both canonically defined as follows (see \cite{BNram} Theorem 3.4):
\begin{lemma}\label{defsec} Let $d, g, r$ be positive integers with  $ r+1-g\le \rho\le 0$.
Choose an admissible  filling of the $(r+1)\times (g-d+r)$-rectangle with the indices $1,2,\dots, g$.
Consider a chain of $g$ elliptic curves $C_1,\dots, C_g$  where $Q_i\in C_i$ glues to $P_{i+1}\in C_{i+1}$ and the curves and points are generic except that when an index $i$ appears twice 
  on spots $(a_1,b_1)$ and $(a_2, b_2)$,  then  $P_i-Q_i$ is a point of torsion precisely  $l_i=|a_2-a_1|+|b_1-b_2|$. 
Then there exist uniquely determined limit sections  $s_a, a=1,\dots, r+1$ that vanish
\begin{itemize}
\item To order $a-1$  at $P_1$ 
\item  To order  $d-ord_{Q_{i-1}}$ at $P_i, i>1$ 
\item To order
$\begin{cases} d-ord_{P_i}s_a&\text{ if } i \text{ is in column a}\\
d-ord_{P_i}s_a-1&\text{ if } i \text{ is not in column a}\end{cases}$\ \ \  at $Q_i$. 
\end{itemize}
The definition for $t_b$ is similar replacing $d$ with $\bar d=2g-2-d$ and $r$ with $\bar r= g-d+r-1$.
\end{lemma}

\begin{lemma}\label{surjP}
With notations as in Lemma \ref{defsec}, choose  $g$ pairs $(a,b)$ such that in the  $\alpha\times \beta$ rectangle the  spots corresponding to the coordinates $(a,b)$  have all different fillings. 
Then, the   $s_a\otimes t_b$ have images in the space of sections of the canonical that are linearly independent.
In particular, if all indices $1,\dots, g$ are used in the filling, the Petri map is surjective.
\end{lemma}
\begin{proof} See  Prop 1.8 in \cite{BNram}
\end{proof}

We will later need the following result:

\begin{lemma}\label{lindep} Let $C$ be an elliptic curve, $P, Q\in C$, $L$ a line bundle of degree $d$ on $C$ 
 and $s_1, \dots, s_{k+1},  k\ge 1$ sections of $L$  that vanish at $P$ to order at least $a$ and at $Q$ with order at least $b$ with $a+b\ge d-k$. 
 Then, $s_1, \dots, s_{k+1}$ are linearly dependent sections. 
\end{lemma}
\begin{proof} By assumption,   $s_1, \dots, s_{k+1}$ give rise to sections of $H^0(C, L(-aP-bQ))$.
By Riemann-Roch's Theorem,  $h^0(C, L(-aP-bQ))=d-a-b\le k$
\end{proof}

\section{Existence of smooth Petri loci of codimension $k$.}

We showed in  \cite{distBN}, Theorem 2.1 that when the Brill Noether number $\rho(g, r, d)=\rho$ is negative but $\rho\ge r+1-g$, 
then there exists smoothable limit linear series on a suitable chain of elliptic curves.
Moreover the Petri map on that chain is onto, so that the kernel of the Petri map has the expected dimension $-\rho$.

We want to show that for $k<\min\{ r+1, g-d+r\}, g\ge (r+1)(g-d+r)-k$, there exist subsets   ${\mathcal P}_{g,d}^{r,k}\subseteq {\mathcal M}_g$ of codimension $k$ and with second Petri map injective.

We will be using the correspondence between limit linear series on a chain of elliptic curves and fillings of an $(r+1)\times (g-d+r)$ rectangle described above.
We assume first that $k=(r+1)(g-d+r)-g=-\rho$.

We define a suitable filling of the $\alpha\times \beta $ rectangle as follows:
We first fill completely the $(\alpha-1)\times (\beta-1)$-rectangle together with the first $\alpha -k-1$ spots of the last row 
and the first $\beta-k-1$ spots of the last column with distinct numbers. 
There are of course many ways to do that, we choose one arbitrarily.
Then we fill the $k$ spots matching terms that appear at the end of the last row with those at the bottom of the last column. 
Finally, we fill the last remaining spot with the index $g$.
As by assumption, $k=(r+1)(g-d+r)-g$, this is possible and the filling uses all the indices $1,\dots, g$.
Moreover,  the indices $i(l)=g-l, l=1,\dots, k$ appear in two spots in the rectangle , namely $(\alpha-l,\beta)$ and $(\alpha,\beta-l)$.
A couple of examples are shown in Figure \ref{fig:puntsdobles}.

\begin{figure}[h!]
\begin{tikzpicture}[scale=.55]

\begin{scope}[xshift=14.5cm]
\foreach \x in {0, 1,2,3,4,5,6,7} {\draw[thick] (0,\x) -- (5, \x); }
\foreach \x in {0, 1,2,3, 4, 5} {\draw[thick] (\x,0) -- ( \x,7); }
\node  [color= blue ] at (.5, 6.5){1}; 
\node[color= blue ] at (1.5, 6.5){7}; 
\node[color= blue ] at (2.5, 6.5){13}; 
\node[color= blue ] at (3.5, 6.5){19}; 
\node[color=blue] at (4.5, 6.5){25}; 

\node[color= blue]  at (.5, 5.5){2}; 
\node[color=blue ] at (1.5, 5.5){8}; 
\node[color=blue ] at (2.5, 5.5){14}; 
\node[color=blue ] at (3.5, 5.5){20}; 
\node[color=blue] at (4.5, 5.5){26}; 

\node [color= blue ]at (.5, 4.5){3}; 
\node[color= blue ] at (1.5, 4.5){9}; 
\node[color= blue ] at (2.5, 4.5){15}; 
\node[color= blue ] at (3.5, 4.5){21}; 
\node[color=purple] at (4.5, 4.5){27}; 

\node [color= blue ] at (.5, 3.5){4}; 
\node  [color=blue] at (1.5, 3.5){10}; 
\node[color= blue ] at (2.5, 3.5){16}; 
\node[color= blue ] at (3.5, 3.5){22}; 
\node[color=purple] at (4.5, 3.5){28}; 

\node [color=blue] at (.5, 2.5){5}; 
\node  [color=blue] at (1.5, 2.5){11}; 
\node[color=blue] at (2.5, 2.5){17}; 
\node[color= blue ] at (3.5, 2.5){23}; 
\node[color=purple] at (4.5, 2.5){29}; 

\node [color=blue]  at (.5, 1.5){6}; 
\node [color=blue] at (1.5, 1.5){12}; 
\node[color=blue] at (2.5, 1.5){18}; 
\node[color=blue] at (3.5, 1.5){24}; 
\node[color= purple ] at (4.5, 1.5){30}; 

\node[color=purple] at (.5, .5){27}; 
\node[color=purple] at (1.5, .5){28}; 
\node  [color=purple] at (2.5, .5){29}; 
\node[color=purple] at (3.5, .5){30}; 
\node[color=blue] at (4.5,.5){31}; 
\end{scope}

\begin{scope}[xshift=20.5cm]
\foreach \x in {0, 1,2,3,4,5,6,7,8,9} {\draw[thick] (0,\x) -- (6, \x); }
\foreach \x in {0, 1,2,3, 4, 5, 6} {\draw[thick] (\x,0) -- ( \x,9); }
\node  [color= blue]  at (.5, 8.5){1}; 
\node[color= blue ] at (1.5, 8.5){9}; 
\node[color= blue ] at (2.5, 8.5){17}; 
\node[color= blue ] at (3.5, 8.5){25}; 
\node[color= blue]  at (4.5, 8.5){33};
\node[color=blue] at (5.5, 8.5){41}; 

\node  [color=blue] at (.5, 7.5){2}; 
\node[color= blue ] at (1.5, 7.5){10}; 
\node[color= blue]  at (2.5, 7.5){18}; 
\node[color= blue ] at (3.5, 7.5){26}; 
\node[color= blue ] at (4.5, 7.5){34}; 
\node[color=blue] at (5.5, 7.5){42}; 

\node [color= blue ]at (.5, 6.5){3}; 
\node[color= blue ] at (1.5, 6.5){11}; 
\node[color= blue]  at (2.5, 6.5){19}; 
\node[color= blue ] at (3.5, 6.5){27}; 
\node[color= blue ] at (4.5, 6.5){35}; 
\node[color=blue] at (5.5, 6.5){43}; 

\node [color= blue ] at (.5, 5.5){4}; 
\node  [color= blue ] at (1.5, 5.5){12}; 
\node[color= blue ] at (2.5, 5.5){20}; 
\node[color= blue ] at (3.5, 5.5){28}; 
\node[color= blue ] at (4.5, 5.5){36}; 
\node[color=blue] at (5.5, 5.5){44}; 

\node [color= blue ] at (.5, 4.5){5}; 
\node  [color= blue ] at (1.5, 4.5){13}; 
\node[color= blue ] at (2.5, 4.5){21}; 
\node[color= blue ] at (3.5, 4.5){29}; 
\node[color= blue ] at (4.5, 4.5){37}; 
\node[color=purple] at (5.5, 4.5){46}; 

\node [color=blue] at (.5, 3.5){6}; 
\node  [color=blue] at (1.5, 3.5){14}; 
\node[color= blue ] at (2.5, 3.5){22}; 
\node[color= blue ] at (3.5, 3.5){30}; 
\node[color= blue ] at (4.5, 3.5){38}; 
\node[color=purple] at (5.5, 3.5){47}; 

\node [color=blue]  at (.5, 2.5){7}; 
\node [color=blue] at (1.5, 2.5){15}; 
\node[color=blue] at (2.5, 2.5){23}; 
\node[color= blue ] at (3.5, 2.5){31}; 
\node[color= blue ] at (4.5, 2.5){39}; 
\node[color=purple] at (5.5, 2.5){48}; 

\node [color=blue]  at (.5, 1.5){8}; 
\node [color=blue] at (1.5, 1.5){16}; 
\node[color=blue] at (2.5, 1.5){24}; 
\node[color=blue] at (3.5, 1.5){32}; 
\node[color= blue ] at (4.5, 1.5){40}; 
\node[color=purple] at (5.5, 1.5){49}; 

\node[color=blue] at (.5, .5){45}; 
\node[color=purple] at (1.5, .5){46}; 
\node  [color=purple] at (2.5, .5){47}; 
\node[color=purple] at (3.5, .5){48}; 
\node[color=purple] at (4.5,.5){49}; 
\node[color=blue] at (5.5, .5){50}; 
\end{scope}

\begin{scope}[xshift=28cm]
\foreach \x in {0, 1,2,3,4,5,6,7} {\draw[thin] (0,\x) -- (5, \x); }
\foreach \x in {0, 1,2,3, 4, 5} {\draw[thin] (\x,0) -- ( \x,7); }
\draw[very thick, purple]  (0,0)--(0,1)--(3,1)--(3,4)--(4,4)--(4,7)--(5,7);

\end{scope}

\end{tikzpicture}
\caption{}
\label{fig:puntsdobles}
\end{figure}

We now look at sections in the kernel of the Petri map.
The section $s_{\alpha -l}\otimes t_{\beta}-s_{\alpha}\otimes t_{\beta -l}$ is in the kernel of the Petri map restricted to $C_{i(l)}$.
We  check the following
\begin{lemma}\label{desker} With the notations in Lemma \ref{defsec} 
 the restriction of the section $s_{\alpha-l}\otimes t_{\beta}-s_{\alpha}\otimes t_{\beta-l}$  in $C_{i( l)}$ glues with sections  in the kernel of the Petri map on each elliptic component 
 giving rise to a section $\bar x_l$ in the kernel of the Petri map on the chain.
\end{lemma}
\begin{proof} From the description of $s_a$ that we gave in Lemma \ref{defsec}  a section $s_a$ vanishes at $P_i$ with order $a-1+i-1-h(i-1,a)$ where $h(i-1,a)$ denotes the number of indices less than $i$ in column $a$.
The section $t_b$ vanishes at $P_i$ with order $b-1+i-1-h'(i-1,b)$ where $h'(i-1,b)$ denotes the number of indices less than  $i$ in row $b$.
Then, the section $s_a\otimes t_b$  vanishes at $P_i$ with order $(a-h'(i-1,b))+(b-h(i-1,a))+2i-2$.

Consider the broken line that limits the filled region of the rectangle with indices up to $i-1$ (the last pane of Figure \ref{fig:puntsdobles} shows this profile 
for the filling in the first pane and $i=22$).
The quantity $(a-h'(i-1,b))+(b-h(i-1,a))$ gives the signed grid distance in the rectangle from the spot $(a,b)$ to this broken line.
In the example of Figure \ref{fig:puntsdobles}, the square $(5, 6)$(last column, last by one row) has horizontal distance 1 and vertical distance 5 to the broken line.
The square $(3,3)$(third row and column) has horizontal distance $-1$ and vertical distance $-3$ to the broken line.

The vanishing at $Q_i$ of the section  $s_a\otimes t_b$ is $2g-2-ord_{P_i}s_a\otimes t_b$ if $i$ is in both column $a$ and row $b$, 
one less if it is only in one of the two and two less if it is in neither.

In particular,  two sections $s_{a_1}\otimes t_{b_1}$ and $s_{a_2}\otimes t_{b_2}$ have the same vanishing at $P_i$
 if the grid distance to the broken line corresponding to $i-1$ is the same. 
They have the same vanishing at $Q_i$ if the grid distance to the broken line corresponding to $i$ is the same. 

Note now that a line bundle on an elliptic curve has at most as many independent sections as its degree, 
except when the degree is 0 when it can have a section if it is the trivial bundle.
Therefore, two sections of a line bundle of degree $2g-2$ with same  orders of vanishing at the points $P_i, Q_i $ that add up to $2g-2$ or $2g-3$ are necessarily linearly dependent,
and the same is true for three sections  when the  orders of vanishing add up to $2g-4$(see Lemma \ref{lindep}).
This is still true if for one of the sections the order of vanishing at $P_i$ or $Q_i$ is one more than for the other sections.

Therefore, to find an element in the kernel of the Petri map restricted to the curve $C_i$, it suffices to find three sections   
with  same distance to the broken lines  corresponding to $i-1$ and $i$. 
It is also enough  to find only two sections $s_a\otimes t_b$ with same distances to the broken lines corresponding to both $i-1$ and $i$
 if in the filling of the rectangle $i$  is in column $a_1$ or in row $b_1$ and also  in column $a_2$ or in row $b_2$. 
Or even if one can find three sections two of them have same distance to the broken line for $i-1$ and distance one less for the line corresponding to $i$.  

We now want to show that the section  $x_l$ which is in the kernel of the Petri map restricted to $C_{i(l)}$ glues with sections in the kernel of the Petri map on  
the curves $C_i, i\not= i(l)$.

On the curves  $C_{i(l')} $ for $  1\le l'\le k$ two squares are added (or removed if $l'>l$) to the last row and last column respectively.
The distances from the squares  $(\alpha-l,\beta )$ to the broken line for the index $i(l')$ is the same as the distance for the square $(\alpha, \beta-l)$. 
Therefore $s_{\alpha-l}\otimes t_{\beta}-s_{\alpha}\otimes t_{\beta-l}$   is in the kernel of the Petri map restricted to  $C_{i( l')}, l'=1,\dots, k$
and these sections for a fixed $l$ glue with each other giving rise to a section $\bar x_l$ in the kernel of the Petri map on the chain for values of $i\ge g-k$.

When going from $C_{i+1}$ to  $C_i$  for  $i_0<i<i(k)$ one square is removed either from the last row or the  last column respectively.
If on the curve $C_{i+1}$, the element in the kernel of the map involved $s_{\alpha}\otimes t_j,  s_j\otimes t_{\alpha}$ 
and the square is removed from the last row, then $s_{\alpha}\otimes t_j, s_{\alpha}\otimes t_{j-1}, s_{j'}\otimes t_{\alpha}$ are linearly dependent.
Similarly, with every new removal, we will obtain a dependence involving two of the sections that appeared on   $C_{i+1}$
and one more pair on the same row or column where the square is being removed.

In going from $C_{i_0+1}$  to $C_{i_0}$, the images of the  sections  $s_{\alpha}\otimes t_1, s_{\alpha-1}\otimes t_{\beta-1}, s_1\otimes t_{\beta}$ are linearly dependent 
and glue with   the extension of $x_{i(l)}$ which was in the kernel of the Petri map for $C_{i_0+1}$.
The  sections  $s_{\alpha}\otimes t_j, s_j\otimes t_{\beta}$  vanish on $P_{i_0}, Q_{i_0}$ to the same order.
So, a linear combination of the two vanishes at $P_{i_0}$ to order one more.
 It follows that a linear combination of $s_{\alpha}\otimes t_j, s_j\otimes t_{\beta}, s_{\alpha}\otimes t_{j-1}, s_{j-1}\otimes t_{\beta}$
is in the kernel of the Petri map and glues with   the extension of $x_{i(1)}$ which was in the kernel of the Petri map for $C_{i_0+1}$.

In going from $C_{i+1}$ to $C_{i}$, if the square we remove in the rectangle did not appear in any of the rows or columns of the sections involved in the element of the kernel of the Petri map in $C_{i+1}$,
then the same sections will be involved in the kernel of the Petri map on $C_{i}$.
If  the square we remove  appears in  the row or column of two of the  sections involved in the element of the kernel of the Petri map in $C_{i+1}$,
then the same sections will be involved in the kernel of the Petri map on $C_{i}$.
If  the square we remove  appears in  one of the rows or columns of two of the  sections involved in the element of the kernel of the Petri map in $C_{i+1}$,
then the same sections will be involved in the kernel of the Petri map on $C_{i}$.
If  the square we remove  appears in  say  the row  of one of the  sections $s_a\otimes t_b$  involved in the element of the kernel of the Petri map in $C_{i+1}$ but not in the row or column of the others, 
then on $C_{i}$. we can use   $s_a\otimes t_b$  and  $s_a\otimes t_{b-1}$ together with the other sections we were using in $C_{i+1}$ to get the element in the kernel of the Petri map. 
\end{proof}

\begin{lemma}\label{difim} With the notations above, the images of 
the sections $\bar x_l,  i=1,\dots, k$  by the second Petri map in the bicanonical linear series are linearly independent. 
\end{lemma}
\begin{proof} Note that on each of the curves $C_{i(1)}=C_{g-k}, \dots, C_{i(k)}=C_{g-1},C_g$, the sections  
\[  s_{\alpha-l}\otimes t_{\beta}-s_{\alpha}\otimes t_{\beta-l},  l=1,\dots, k\]
  are in the kernel of the Petri map.
It suffices to check that on   $C_g$ the images of these sections on the bicanonical are linearly independent.
 
On $C_g$, the section $s_a$ vanishes at $P_g$ to order $a+g-\beta-2$  if $a\not=\alpha$ and one more if  $a= \alpha$ and at $Q_g$ to order $d-a-g+\beta +1$.  

Therefore, the two sections $s_{\alpha -l},  s_{\alpha}$ span a linear series with fixed part 
\[ (\alpha-l-\beta-1)P_g+(d-\alpha-g+\beta+1)Q_g=(d-1-l)P_g\]
 (where we used that   $\alpha=r+1, \beta=g-d+r$).
 The mobile part is spanned by two sections of   $(l+1)P_g$ one with orders of vanishing $l+1$ at $P_g$  and 0 at $Q_g$ and the other with 
  order of vanishing $l$ at $Q_g$  and 0 at $P_g$
 Writing $e$ for a section of the fixed part, then  $s_{\alpha}=es'_{\alpha}, s_{\alpha -l}=es'_{\alpha -l}$. Write $f$ for the   function such that $s'_{\alpha -l}=fs'_{\alpha}$.

Similarly, a section $t_a$ vanishes at $P_{g}$ to order $a+g-\alpha-2$  if $a\not=\beta$ and one more if  $a= \beta$ and at $Q_g$ to order $\bar d-a-g+\alpha +1 $ where $\bar d=2g-2-d$.
Therefore, the two sections $t_{\beta -l},  t_{\beta}$ span a linear series with fixed part  $(\bar d-1-l)P_g$
and the mobile part is spanned by two sections of   $(l+1)P_g$ one with orders of vanishing $l+1$ at $P_g$  and 0 at $Q_g$ and the other with 
  order of vanishing $l$ at $Q_g$  and 0 at $P_g$ (as was to be expected from the base-point  free pencil trick).
 Writing $\bar e$ for a section of the fixed part, then  $t_{\beta}=\bar et'_{\beta}, t_{\beta -l}=\bar et'_{\beta -l}$. 
 Then also   $t'_{\beta -l}=ft'_{\beta}$ for the  same  function $f$ as before.

The image of $s_{\alpha -l}\otimes t_{\beta}-s_{\alpha}\otimes t_{\beta -l}=efs'_{\alpha}\otimes \bar e t'_{\beta}-es'_{\alpha}\otimes \bar e ft'_{\beta}$ in the bicanonical series  
is the section with  divisor $[(d-1-l)+(\bar d-1-l)]P_g+R_f$ with $R_f$ the ramification divisor of the function $C\to {\mathbb P}^1$ defined by $f$ 
(see the argument for a $g^1_k$ at the end of section \ref{secdeff} or Lemma 2.12 of \cite{semican}).
 In this case, the ramification divisor of the moving portion is $lP_g+(l-1)Q_g+R'$ with $R'$ a divisor not supported at $P_g, Q_g$.
 Therefore, the image in the bicanonical is a section with divisor $(2g-4-l)P_g+lQ_g+R$.
 For different values of $l$, these sections have different orders of vanishing at   $P_g,Q_g$, therefore they are linearly independent.
\end{proof}

\begin{lemma}\label{anulQ} With the notations above,
the sections $\bar x_l,  l=1,\dots, k$ vanish to order $l$ at $Q_g$. 
\end{lemma}
\begin{proof}  Proved above.
\end{proof}

We now use the above results to construct limit sections in the kernel of the Petri map for curves with higher Brill-Noether number.
The construction is similar to  the one used in the proof of Theorem 1.3 in  \cite{F2} that showed the existence of Petri divisors with an inductive method. 

\begin{lemma}\label{ggeneral} Let $g, r, d, k$ be integers such that $-\rho\le k\le r, r+1\le g-d+r$.
There exists then a family of codimension $k$ in the space of  chains of elliptic curves with a limit linear series on it which has kernel of the Petri map precisely $k$.
\end{lemma}
\begin{proof}  Define $t=\rho+k, g_0=g-t, d_0=d-t$.
Then $\rho(g_0, d_0, r)=\rho(g, d, r)-t=-k$.
From Lemma \ref{difim} , there exists a subset of codimension $k$ of the chains of $g_0$ elliptic curves 
with limit linear series $(L_i, V_i), i=1,\dots, g_0$ of degree $d-t$ and dimension $r$ with kernel of the Petri map of dimension $k$.

Consider one such chain and limit linear series.
Add to to the chain  $t$ generic elliptic curves $C_{g_0+1},\dots, C_g$ with generic points $P_i, Q_i\in C_i, i=g_0+1,\dots, g$ by attaching $ Q_{i}\in C_{i} $ to $ P_{i+1}\in C_{i+1}, i=g_0,\dots, g-1 $.

On the curves $C_i, i=1,\dots, g_0$ take the limit linear series with aspects $(L_i(tQ_i), V_i+tQ_i)$.
On each of the curves $C_i, i=g_0+1,\dots, g$ choose a line bundle ${L'_i}$ of degree $r+1$ such that ${L'_i}^2={\mathcal O}(2(r+1)P_i)$ but ${L'_i}\not={\mathcal O}((r+1)P_i)$.
Take 
\[  L_i={L'_i}((d-r-g+i-1)P_i+(g-i)Q_i), \ V_i=|{L'}_i|+(d-r-g+i-1)P_i+(g-i)Q_i.\]

The Serre dual linear series on the curves  $C_i, i=1,\dots, g_0$ is given by  $(\overline{L_i}(tP_i), \overline{V_i}+tP_i)$ 
where    $(\overline{L_i}, \overline{V_i})$ were the Serre dual of the original linear series in the original chain of $g_0$ curves..
On $C_i, i=g_0+1,\dots, g$ the Serre dual bundle and  space of section  are 
\[ L'_i((g-d+i-r-3)P_i+(g-i)Q_i),\ W_i  \text{ such that } |L'_i|+(g-d+i-r-3)P_i+(g-i)Q_i\subseteq W_i\]

From Lemma \ref{anulQ}, the restriction of the limit linear series to the first $g_0$ components has $k$ sections in the kernel of the Petri map 
vanishing at $Q_{g_0}$ to orders $l+2(g-g_0)=l+2t,\ l=0,\dots, k$.

The linear series $|{L'}_i|$ has sections $s_j, j=0,\dots, r$ vanishing at $P_i$ to order $r-j$ and at $Q_i$ to order $j$.
If $e$ denotes a section vanishing to order $d-r-g+i-1$ at $P_i$ and $g-i$ at $Q_i$ then $es_j$ is a section of the linear series $|{L'}_i|+(d-r-g+i-1)P_i+(g-i)Q_i$.
Similarly, if  $\bar e$ denotes a section vanishing to order $g-d+i-r-3$ at $P_i$ and $g-i$ at $Q_i$ then $\bar es_j$ is a section of the linear series $ |L'_i|+(g-d+i-r-3)P_i+(g-i)Q_i\subseteq W_i $.
 
 The section $es_r\otimes \bar e s_{r-l}-es_{r-l}\otimes \bar e s_{r }$ on the curve $C_i $ is in the kernel of the Petri map and vanishes at $P_i$ to order 
 \[ (d-r-g+i-1)+(  (g-d+i-r-3)+r+(r-l)=2i-4-l .\]
 The section $es_r\otimes \bar e s_{r-l}-es_{r-l}\otimes \bar e s_{r }$ on the curve $C_{i-1} $ vanishes at $Q_{i-1}$ to order 
 \[ (g-(i-1))+(g-(i-1))+l=2g-2i+2+l. \]
As $(2i-4-l )+(2g-2i+2+l)=2g-2=\deg K_C$,  the sections $es_j\otimes \bar e s_{j'}-es_{j'}\otimes \bar e s_{j}$ can be glued to form a section $\hat x_l$ in the chain $g_0+1\le i\le g$.
   From Lemma \ref{anulQ}     the section obtained from $\bar x_l$ by adding $tQ_i, i=1,\dots, g_0$ to each of its aspects can be glued with $\hat x_l$ to give a section $\tilde x_l$ in the kernel of the Petri map on the whole chain
  By construction, the restrictions  to $C_{g_0}$ of $\tilde x_l$ is  $\bar x_l +tQ_{g_0}$, so from Lemma \ref{difim}, their images in the bicanonical are linearly independent.
\end{proof}

\begin{theorem}\label{ggeneral} Let $g, r, d, k$ be integers such that $-\rho\le k\le r, r+1\le g-d+r$.
There exists then a component of ${\mathcal P}_{g,d}^{r,k}$ of codimension $k$ in ${\mathcal M}_g$.
\end{theorem}
\begin{proof}  Let $C_o$ be the chain of elliptic curves constructed in Lemma \ref{ggeneral} .
Denote by $\pi: {\mathcal C}\to B$ the versal deformation space of $C_0$ with $b_0\in B$ corresponding to $b_0$.
From Theorem 1.2 of \cite{F2}, there exists a quasi-projective variety  ${\mathcal Q}_{g,d}^{r,1}$ parameterizing limit linear series on the fibers of $\pi$ together with a section of these limit section in the kernel of the Petri map.
Moreover, every irreducible component of ${\mathcal Q}_{g,d}^{r,1}$ has dimension at least $\dim B-1$.
By construction, the fiber of $p:{\mathcal Q}_{g,d}^{r,1}\to B$ over $b_0$ has dimension $k$ with $(C_0, g ^r_d,\tilde x_l)$ spanning this fiber.
The tangent space at $b_0$ to the image of $p$ is the orthogonal to the image by  $\mu_1$ of the images of the span  of the $\tilde x_l$ which has codimension $k$.
Therefore, the space where all $k$ sections deform has codimension $k$ in ${\mathcal M}_g$.

\end{proof}



\begin{thebibliography}{ACGH85}

\bibitem[AC]{AC} E. Arbarello,, M. Cornalba, 
{\it Su una congettura di Petri}  Comment. Math. Helv.56(1981), no.1, 1-38.

\bibitem[B]{Benzo} L.~Benzo, 
{\it Components of moduli spaces of spin curves with the expected codimension}, Math. Ann. {\bf 363} (2015), no.~1-2, 385--392. 


\bibitem[BS]{BS} A. Bruno, E. Sernesi
{\it A note on the Petri loci} Man. Math. 136, (2011)  439-443.

\bibitem[C]{C} A.Castorena,  
{\it Remarks on the Gieseker-Petri divisor in genus eight}
Rend. Circ. Mat. Palermo (2) 59 (2010), no. 1, 143-150.

\bibitem[CT]{CT} A.Castorena,  M. Teixidor, 
{\it  Divisorial components of the Petri locus for pencils}
J. Pure Appl. Algebra 212 (2008), no. 6, 1500–1508.

\bibitem[E]{E}
D. Edidin, {\it The monodromy of certain families of linear series is at least the alternating group.},
 Proc. Amer. Math. Soc.  \textbf{113} (1991), no.~4, 911--922.
 
 \bibitem[EH]{EH} D.~Eisenbud, J.~Harris, {\it Irreducibility of some families of linear series with Brill-Noether number -1}, Ann. Sci. École Norm. Sup. (4)22(1989), no.1, 33-53.
 
\bibitem[F1]{F1} G. Farkas.{\it  Gaussian maps, Gieseker–Petri loci and large theta-characteristics.} J. Reine
Angew. Math. 581, 151-173 (2005)

\bibitem[F2]{F2} G. Farkas.{\it   Rational maps between moduli spaces of curves and Gieseker–Petri divisors.} J. Algebr. Geom. 19  (2010), 243-284.


\bibitem[HT]{HT} R. Haburcak, M.~Teixidor, {\it Some reducible and irreducible Brill-Noether loci}
arXiv:2503.16255  

 \bibitem[H]{H}
J. Harris,  {\it Theta-characteristics on algebraic curves}.
Trans. Amer. Math. Soc. 271 (1982), no. 2, 611–638.
 


\bibitem[LC1]{LC1} M. Lelli-Chiesa,
{\it  The Gieseker-Petri divisor in  ${\mathcal M}_g$ for  $g\le 13$} 
Geom. Dedicata 158 (2012), 149-165.

 \bibitem[LC2]{LC2} M. Lelli-Chiesa,
{\it  A codimension 2 component of the Gieseker-Petri locus.}
J. Algebraic Geom. 31 (2022), no. 4, 751–771.

\bibitem[T1]{semican} M.~Teixidor, {\it Half-canonical series on algebraic curves}, Trans. Amer. Math. Soc.302, n1 (1987), 9-115.

\bibitem[T2]{BNram}
 M.~Teixidor.
 {\it Brill-Noether loci with ramification at two points}. Ann Mat. Pura Appl {\bf 202} 2023 n3, 1217-1232.


\bibitem[T3]{distBN}
 M.~Teixidor.
\newblock{ \it Brill-Noether Loci},
\newblock Manuscripta Math. 176 (2025), no. 1,In press. 











\end{thebibliography}
\end{document}